\documentclass[leqno,12pt]{article}

\usepackage{amssymb,amsmath,amsfonts,latexsym}
\newtheorem{theorem}{Theorem}

\newtheorem{lemma}{Lemma}

\def\Z{\ensuremath{\mathbb{Z}}}
\def\Q{\ensuremath{\mathbb{Q}}} 
\def\P{\ensuremath{\mathbb{P}}}

\def\int{\mathrm{int}}

\def\conv{\mathrm{conv}}
\def\Pic{\ensuremath{\mathrm{Pic}}}

\def\<{\ensuremath{\langle}}
\def\>{\ensuremath{\rangle}}

\title{Fujita's very ampleness conjecture for singular toric varieties}

\author{\textsc{Sam Payne}}

\date{}

\begin{document}

\maketitle

\begin{abstract}
    
We present a self-contained combinatorial approach to Fujita's
conjectures in the toric case.  Our main new result is a
generalization of Fujita's very ampleness conjecture for toric
varieties with arbitrary singularities.  In an appendix, we use
similar methods to give a new proof of an analogous toric
generalization of Fujita's freeness conjecture due to Fujino.

\end{abstract}

\section{Introduction}

Given an ample divisor $D$ and any other Cartier divisor $D'$ on an
algebraic variety, we can choose $t$ sufficiently large so that $tD +
D'$ is basepoint free or very ample.  In either case, it is not easy
to say how large we must choose $t$ in general.  However, for the case
where $D'$ is in the canonical class $K_{X}$, Fujita made the
following conjectures.

\vspace{10 pt}

\noindent \textbf{Fujita's Conjectures} \hspace{0 pt} \emph{Let $X$ be
an $n$-dimensional projective algebraic variety, smooth or with mild
singularities, $D$ an ample divisor on $X$.}

     (i) \emph{For $t \geq n + 1,\ tD + K_{X}$ is basepoint free.}

     (ii) \emph{For $t \geq n + 2,\ tD + K_{X}$ is very ample.}

\vspace{10 pt}     
     
\noindent The case where $X$ is $\P^{n}$ and $D$ is a hyperplane shows
that Fujita's conjectured bounds are best possible.

For smooth varieties, the corresponding statements with ``basepoint
free'' and ``very ample'' replaced by ``nef'' and ``ample'',
respectively, are consequences of Mori's Cone Theorem \cite{Fujita}. 
For divisors on smooth \emph{toric} varieties, nefness and ampleness
are equivalent to freeness and very ampleness, respectively, so
Fujita's conjectures follow immediately for smooth toric varieties. 
One can also deduce Fujita's conjectures for smooth toric varieties by
general (non-toric) cohomological arguments of Ein and Lazarsfeld in
characteristic zero \cite{Ein-Lazarsfeld}, and Smith in positive
characteristic \cite{Smith1} \cite{Smith2}, again using the fact that
ample divisors on smooth toric varieties are very ample.

For toric varieties with arbitrary singularities, a strong
generalization of Fujita's freeness conjecture was proved by
Fujino \cite{Fujino}.  We follow the usual toric convention fixing
$K_{X} = - \sum D_{i},$ the sum of the $T$-invariant prime divisors
each with coefficient -1, as a convenient representative of the
canonical class.

\vspace{10 pt}

\noindent \textbf{Fujino's Theorem} \hspace{0 pt} \emph{Let $X$ be a
projective $n$-dimensional toric variety not isomorphic to $\P^{n}$. 
Let $D$ and $D'$ be \Q -Cartier divisors such that $0 \geq D' \geq K_{X}$,
$D+D'$ is Cartier, and $D \cdot C \geq n$ for all $T$-curves $C$. 
Then $D + D'$ is basepoint free.}

\vspace{10 pt}

Fujita's freeness conjecture for toric varieties is the
special case of Fujino's Theorem when $X$ is Gorenstein, $D' = K_{X}$,
and $D = t L$ for some ample Cartier divisor $L$ and some integer $t
\geq n + 1$.  Fujino's Theorem shows that, for toric
varieties, Fujita's conjectured bound can be improved by
excluding the extremal case $X \cong \P^{n}$.  

Of course, the case of $\P^{n}$ can be analyzed separately.  The
canonical divisor on $\P^{n}$ is linearly equivalent to $(-n-1)H$, so
$D' \sim sH$ for some $0 \geq s \geq (-n-1).$ Any \Q -Cartier divisor
$D$ on $\P^{n}$ is linearly equivalent to $tH$ for some $t \in \Q$. 
Then $D+D'$ is Cartier exactly when $t + s$ is an integer and
basepoint free exactly when $t + s$ is a nonnegative integer.

The main purpose of this paper is to prove an analogous
generalization of Fujita's very ampleness conjecture for toric
varieties with arbitrary singularities.

\begin{theorem} \label{Main} Let $X$ be a projective $n$-dimensional
toric variety not isomorphic to $\P^{n}$.  Let $D$ and $D'$ be \Q -Cartier
divisors such that $0 \geq D' \geq K_{X}$, $D+D'$ is Cartier, and $D
\cdot C \geq n + 1$ for all $T$-curves $C$.  Then $D + D'$ is very
ample.
\end{theorem}

The statement of Fujino's Theorem can be strengthened by removing the
assumption that $D+D'$ is Cartier.  A Cartier divisor on a toric
variety is basepoint free if and only if it is nef, i.e., if and only
if $D$ intersects every curve nonnegatively \cite[Proposition
1.5]{Laterveer}.  Without the hypothesis that $D+D'$ is Cartier, the
sharper statement of Fujino's Theorem, which one may deduce from
\cite[Theorem 0.1]{Fujino}, is then:

\vspace{10 pt}

\noindent \textbf{Fujino's Theorem$^{+}$} \hspace{0 pt}
\emph{Let $X$ be a projective $n$-dimensional toric variety not
isomorphic to $\P^{n}$.  Let $D,D'$ be \Q -Cartier divisors such
that $0 \geq D' \geq K_{X}$ and $D \cdot C \geq n$ for all $T$-curves
$C$.  Then $D +D'$ is nef.}

\vspace{10 pt}

Similarly, the statement of Theorem 1 can be strengthened by using a
toric characterization of very ampleness to remove the hypothesis that
$D+D'$ is Cartier.  Since every divisor on a toric variety is linearly
equivalent to a $T$-invariant divisor, we may assume $D$ and $D'$ are
\Q -linear combinations of $T$-invariant divisors.  To state the 
stronger theorem, we need some notation from toric geometry.  

Let $D = \sum d_{i}D_{i}$ be a $T$-\Q -Cartier divisor on a complete
toric variety $X$.  Let $M$ the character lattice of $T$, and let
$M_{\Q} = M \otimes \Q$.  For each maximal cone $\sigma$ in the fan
defining $X$, we have a point $u_{\sigma} \in M_{\Q}$ determined by
the conditions $\< u_{\sigma}, v_{i} \> = -d_{i}$ for each of the
primitive generators $v_{i}$ of the rays of $\sigma$.  When $D$ and
$D'$ denote $T$-\Q -Cartier divisors, we will write $u_{\sigma}$ and
$u'_{\sigma}$ for the points of $M_{\Q}$ associated to $D$ and $D'$,
respectively.  The association $D \leadsto u_{\sigma}$ is linear,
i.e., $tD \leadsto t u_{\sigma} $ and $D+D' \leadsto u_{\sigma} +
u'_{\sigma} $.  A $T$-\Q -Cartier divisor $D$ is Cartier if and only
if $u_{\sigma} \in M$ for all maximal cones $\sigma$.  Also associated
to $D$ is a polytope $P_{D} \subset M_{\Q}$ cut out by the
inequalities $\<u, v_{i} \> \geq -d_{i}$ for all of the primitive
generators $v_{i}$ of the rays of the fan.  When $D$ is $T$-\Q
-Cartier and nef, the $ \{ u_{\sigma} \} $ are the vertices of
$P_{D}$.  If we translate $P_{D}$ so that the vertex $u_\sigma$ is at
the origin, then all of $P_{D}$ sits inside the dual cone
$\sigma^{\vee}$.  Write $P_{D}^{\sigma}$ for this translation, i.e.,
$P_{D}^{\sigma} := P_{D} - u_{\sigma}$.  In case $D$ is Cartier, then
$D$ is very ample if and only if $P_{D}^{\sigma} \cap M$ generates the
semigroup $\sigma^{\vee} \cap M$ for all maximal cones $\sigma$.
Without the hypothesis that $D+D'$ is Cartier, the stronger version of
Theorem \ref{Main} we will prove is then:

\begin{theorem}\label{Main'} Let $X$ be a projective $n$-dimensional
toric variety not isomorphic to $\P^{n}$.  Let $D$ and $D'$ be $T$-\Q -Cartier
divisors such that $0 \geq D' \geq K_{X}$ and $D \cdot C \geq n + 1$
for all $T$-curves $C$.  Then $P_{D+D'}^{\sigma} \cap M$ generates
$\sigma^{\vee} \cap M$ for all maximal cones $\sigma$.
\end{theorem}

Theorem 1 is the special case of Theorem 2 when $D+D'$ is 
Cartier.

\vspace{5 pt}

Our approach starts from an observation made by Laterveer in
\cite{Laterveer}: if $D$ is ample, then the lattice length of the edge
of $P_{D}$ corresponding to a $T$-curve $C$ is precisely $D \cdot
C$.  This fact can be seen as a consequence of Riemann-Roch
for toric varieties \cite[p.112]{Fulton}.  Adding $D'$ corresponds to
moving the faces of $P_{D}$ inward at most a unit distance with
respect to the dual lattice.  When all of the edges of $P_{D}$ have
lattice length at least $n + 1$, we show that $P_{D+D'}^{\sigma}$
contains an explicit generating set for $\sigma^{\vee} \cap M$ for all
maximal cones $\sigma$.  The computations are straightforward in the
simplicial case, as can be seen in the example at the end of the
introduction.

In the proof of Fujino's Theorem, there is a simple reduction to the
simplicial case (see \cite[Lemma 2.4]{Laterveer} or \cite[1.12, Step
2]{Fujino}).  That reduction works via a partial projective resolution
of singularities corresponding to a regular triangulation of the fan
defining $X$.  This type of reduction seems not to work for very
ampleness.  Instead, for each nonsimplicial maximal cone $\sigma$, we
make a canonical subdivision of the dual cone $\sigma^{\vee}$.

\vspace{5 pt}

The earliest results on Fujita's freeness conjecture for singular
toric varieties of which the author is aware are due to Laterveer.  In
\cite{Laterveer}, Laterveer proved Fujino's Theorem for \Q -Gorenstein
toric varieties when $D' = K_{X}$ using toric Mori theory, as
developed in \cite{Reid}.  Our statement of Theorem \ref{Main}, like
the statement of Fujino's Theorem, is influenced by Musta\c{t}\v{a}'s
formulations in \cite{Mustata}.  In particular, Musta\c{t}\v{a} stated
and proved Fujino's Theorem and Theorem \ref{Main} for smooth toric
varieties when $D$ and $D'$ are Cartier as consequences of a
characteristic-free vanishing theorem for toric varieties.  For proofs
of Fujino's Theorem that do not use vanishing theorems or toric Mori
theory, see also \cite{Lin} or the Appendix.

The only previous result on Fujita's very ampleness conjecture for
singular toric varieties of which the author is aware is due to
Lin,\footnote{Another effective very ampleness result for singular
toric varieties, due to Ewald and Wessels \cite{Ewald}, may be stated
as follows: let $X$ be an $n$-dimensional projective toric variety and
$D$ a $T$-\Q -Cartier divisor on $X$ such that $D \cdot C \geq n-1$
for all $T$-curves $C$.  Then $P^{\sigma}_{D} \cap M$ generates
$\sigma^{\vee} \cap M$ for all maximal cones $\sigma$.  In particular,
if $D$ is Cartier, then $D$ is very ample.} who proved the conjecture
for simplicial Gorenstein toric varieties in dimension $\leq 6$
\cite{Lin}.  In \cite{Laterveer}, Laterveer also claimed to prove a
generalization of Fujita's very ampleness conjecture for arbitrary \Q
-Gorenstein toric varieties.  As noted by Lin, there is an error in
the proof of this claim.  In particular, it is not true in general
that $P_{tD + K_{X}}$ contains $P_{(t-1)D}$.  The case where $X$ is
$\P^{n}$ and $D$ is a hyperplane is a counterexample.  Nevertheless,
Laterveer's approach to Fujita's very ampleness conjecture for
singular toric varieties contains fruitful insights, in particular,
the realization that $\P^{n}$ is the only toric extremal case and the
characterization of the intersection numbers $D \cdot C$, for
$T$-curves $C$, as the lattice lengths of the edges of $P_{D}$.  The
results we prove here are strong enough to imply all of the very
ampleness results claimed in \cite{Laterveer}.

\vspace{5 pt}

\noindent \textbf{Example} We illustrate the essential techniques of
this paper in a concrete simplicial example.  Let $u_{1}, \ldots,
u_{n}$ be linearly independent primitive vectors in $M$.  Let $P$ be
the simplex with vertices $\{0, u_{1}, \ldots , u_{n} \}$.  Associated
to $P$, there is a projective toric variety $X_{P}$ with an ample
divisor $D$ such that $P = P_{D}$ \cite[Section 1.5]{Fulton}.  The
vertex $0$ of $P$ corresponds to a maximal cone $\sigma$ of the fan
defining $X_{P}$ whose dual cone $\sigma^{\vee}$ is spanned by $\{u_1,
\ldots , u_{n}\}$.  Let $D_{1},\ldots, D_{n}$ be the divisors
corresponding to the rays of $\sigma$, and let $D' = -D_{1}- \cdots
-D_{n}$.  We will show that, for $t \geq n + 1$, $P_{tD+D'}^{\sigma}
\cap M$ generates $\sigma^{\vee} \cap M$.

Every point in $M$ can be written uniquely as an integer linear
combination of the $\{ u_{i} \}$ plus a fractional part.  So the
semigroup $\sigma^{\vee} \cap M$ is generated by $\{ 0, u_{1}, \ldots,
u_{n} \}$ together with $\{ (a_{1}u_{1} + \cdots + a_{n} u_{n}) \in M
\ | \ 0 \leq a_{i} < 1 \}$.  For $t \geq n + 1$, we will show that
$P_{tD+D'}^{\sigma}$ contains this generating set.

Define a linear function $\lambda$ on $M_{\Q}$ by \[
\lambda(a_{1}u_{1} + \cdots + a_{n}u_{n}) \ = \ a_{1} + \cdots +
a_{n}.  \] Note that $P_{tD} = \{ u \in \sigma^{\vee} | \lambda(u)
\leq t \}.$ In other words, if $\{ v_{i} \}$ are the primitive
generators of the rays of $\sigma$, then $P_{tD}$ is cut out by the
conditions $\< u , v_{i} \> \geq 0$ and the condition $\lambda(u) \leq
t$.  Similarly, $P_{tD+D'} = \{u'_{\sigma} + u \ | \ u \in
\sigma^{\vee}, \lambda( u'_{\sigma} + u ) \leq t \}$, i.e.\
$P_{tD+D'}$ is cut out by the conditions $\< u, v_{i} \> \geq \<
u'_{\sigma}, v_{i} \> = 1$ and the condition $\lambda(u) \leq t$.  It
follows that any lattice point in $P_{tD}$ that is in the interior of
$\sigma^{\vee}$ is contained in $P_{tD+D'}$.  Indeed, if $u$ is in
$P_{tD}$, then $\lambda(u) \leq t$, and if $u$ is a lattice point in
the interior of $\sigma^{\vee}$, then $\< u , v_{i} \>$ is a positive
integer.

Suppose $t \geq n + 1$.  Then $u_{1} + \cdots + u_{n}$ is a lattice 
point in $P_{tD}$ that is in the interior of $\sigma^{\vee}$, so 
$u_{1} + \cdots + u_{n}$ is contained in $P_{tD + D'}$.  Note that 
$u'_{\sigma}$ is the point of $P_{tD+D'}$ for which $\lambda$ 
achieves its minimum.  In particular, $\lambda(u'_{\sigma}) \leq 
\lambda(u_{1} + \cdots + u_{n}) = n$.

For each $u_{i}$, we have $\lambda(u'_{\sigma} + u_{i}) =
\lambda(u'_{\sigma}) + 1 \leq n + 1 \leq t$.  Therefore $u'_{\sigma} +
u_{i} \in P_{tD+D'}$, i.e., $u_{i} \in P_{tD+D'}^{\sigma}$.  Given a
lattice point $p$ of the form $p = a_{1} u_{1} + \cdots + a_{n}u_{n}$
with $0 \leq a_{i} < 1$, we have another lattice point $p' = (1-a_{1})
u_{1} + \cdots + (1-a_{n})u_{n}$ in $P_{tD}$ that is in the interior
of $\sigma^{\vee}$.  So $p'$ is contained in $P_{tD + D'}$, and
therefore \[ \lambda(u_{\sigma}') \ \leq \ \lambda(p') \ = \ n - \lambda 
(p). \]
So $\lambda(u_{\sigma}' + p) \leq n < t$, and hence $p \in
P_{tD+D'}^{\sigma}$, as required.

\vspace{5 pt}

I wish to thank M. Hering, P. Horja, R. Lazarsfeld, and M.
Musta\c{t}\v{a} for helpful conversations related to this work.  I am
especially grateful to W. Fulton for his encouragement on this project
and for his comments and suggestions on earlier drafts of this paper.

\section{Preliminaries}

As a first step to proving Theorem 2, we have:

\begin{lemma} \label{cor to Fuj} Let $X, D$, 
and $D'$ satisfy the hypotheses of Theorem {\rm 2}.  Let $\sigma$ be a 
maximal cone, and let $\{ u_{1}, \ldots , u_{s} \}$ be the primitive 
generators of the rays of $\sigma^{\vee}$.  Then
$P^{\sigma}_{D+D'}$ contains $\{0, u_{1}, \ldots, u_{s}\}.$ 
\end{lemma}    

Proof. By Fujino's Theorem$^{+}$, $(n/(n+1))D + D'$ is nef. 
Therefore, for any $T$-curve $C$, \[ (D+D') \cdot C \ = \ \frac{1}{n+1} (D
\cdot C) + (\frac{n}{n+1}D + D') \cdot C \ \geq \ 1.\] By Laterveer's
observation, this means that every edge of $P_{D+D'}$ has lattice
length at least 1.  Translating the vertex $u_{\sigma}$ to the origin,
it follows that $P_{D+D'}^{\sigma}$ contains 0 and the primitive
generators of each of the rays of $\sigma^{\vee}$. \hfill $\Box$

\vspace{5 pt}

If $\sigma$ is regular, then $s = n$ and $\{ 0, u_{1}, \ldots, u_{n}
\}$ generates $\sigma^{\vee} \cap M$, so the conclusion of Theorem 2,
i.e., the fact that $P^{\sigma}_{D+D'} \cap M$ generates $\sigma^{\vee}
\cap M$, follows immediately.  In general, if we let $\Delta = \conv
\{ 0, u_{1}, \ldots, u_{s} \}$, then Lemma 1 says that
$P_{D+D'}^{\sigma}$ contains $\Delta$.  If $\sigma$ is not regular,
then $\Delta$ may not contain a generating set for $\sigma^{\vee} \cap
M$.  The following example, due to Ewald and Wessels \cite{Ewald},
illustrates this possibility.

\vspace{5 pt}

\noindent \textbf{Example} Let $M = \Z^{3}; u_{1} = (1,0,0), u_{2} =
(0,1,0),$ and $u_{3} = (1,1,2)$.  Let $\sigma^{\vee}$ be the cone
spanned by $\{ u_{1}, u_{2}, u_{3} \}$, so $\Delta = \conv \{ 0,
u_{1}, u_{2}, u_{3} \}$.  Then $\Delta \cap M = \{ 0 , u_{1} , u_{2},
u_{3} \}$, so the semigroup generated by $\Delta \cap M$ only contains
lattice points whose third coordinate is even.  In particular, the
lattice point $(1,1,1) = (1/2)(u_{1} + u_{2} + u_{3})$ is in
$\sigma^{\vee}$, but not in the semigroup generated by $\Delta \cap
M$.

\vspace{5 pt}

Although $\Delta$ may not contain a generating set for $\sigma^{\vee}
\cap M$, we will show that $P_{D+D'}^{\sigma}$ contains a dilation of
$\Delta$ that does contain a generating set.  Let $m = \min \{ (D+D')
\cdot V(\sigma \cap \tau) \}$, where $\tau$ varies over all maximal 
cones adjacent to $\sigma$, so that $m$ is the minimum of the
lattice lengths of the edges of $P^{\sigma}_{D+D'}$ incident to the
vertex 0.  Note that $m \Delta$ is the largest rational dilation of
$\Delta$ contained in $P_{D+D'}^{\sigma}$.  We will show that $m
\Delta$ does contain a generating set for $\sigma^{\vee} \cap M$.

In preparation for proving this, we develop a few preliminaries.  
First, we generalize Laterveer's observation on the lattice
lengths of the edges of $P_{D}$ to the case where $D$ is not 
necessarily ample.

\begin{lemma}\label{Laterveer} Let $X$ be a complete toric variety and
$D$ a $T$-\Q -Cartier divisor on $X$.  Let $\sigma, \tau$ be adjacent
maximal cones in the fan defining $X$, and let $u$ be the primitive
generator of the ray of $\sigma^{\vee}$ perpendicular to $\sigma \cap
\tau$.  Then \[ u_{\tau} \, = \, u_{\sigma} + (D \cdot V(\sigma \cap \tau))
u.  \]
\end{lemma}

Proof. Since $u_{\sigma}$ and $u_{\tau}$ agree on $\sigma \cap \tau$,
their difference must vanish on $\sigma \cap \tau$, i.e.\ $u_{\tau} -
u_{\sigma} = k u$ for some rational number $k$.  By the toric
intersection formulas in \cite[Section 5.1]{Fulton}, for $v_{j}$ the
primitive generator of any ray of $\tau$ not contained in $\sigma$, \[
D \cdot V(\sigma \cap \tau) \ = \ \frac{ \< u_{\sigma} - u_{\tau} ,
v_{j} \>}{- \< u , v_{j} \>}.  \] Therefore, \[ D \cdot V(\sigma \cap
\tau) \ = \ \frac{ \< ku, v_{j} \> } { \< u, v_{j} \> } = k.  \]

\vspace{-20 pt} \hfill $\Box$

\vspace{30 pt}

Now we develop some tools for working with rational cones.  Let
$\sigma^{\vee}$ be a strictly convex $n$-dimensional rational cone,
and let $u_{1}, \ldots, u_{s}$ be the primitive generators of the rays
of $\sigma^{\vee}$.  Define a function $\lambda^{\min}$ on
$\sigma^{\vee}$ by \[ \lambda^{\min}(u) \ = \ \min \{ (a_{1} + \cdots
+ a_{s}) \ | \ a_{1} u_{1} + \cdots + a_{s} u_{s} = u, a_{i} \geq 0 \}.
\] Define $\lambda^{\max}$ similarly.  A few combinatorial properties
of $\lambda^{\min}$ and $\lambda^{\max}$, all of which are immediate from
the definitions, will be useful in what follows.

First, $\lambda^{\min}$ and $\lambda^{\max}$ are anticonvex
and convex, respectively.  In other words, for any $u, u' \in
\sigma^{\vee}$, \[ \lambda^{\min}(u+u') \ \leq \ \lambda^{\min}(u) +
\lambda^{\min}(u'), \] and similarly $\lambda^{\max}(u+u') \geq
\lambda^{\max}(u) + \lambda^{\max}(u')$.  

Second, suppose the restriction of $D' = \sum d'_{i}D_{i}$ to the
affine open $U_{\sigma}$ is minus-effective, i.e.\ for each of the
primitive generators $v_{i}$ of the rays of $\sigma$, $d'_{i} \leq 0$. 
Then $\< u'_{\sigma}, v_{i} \> = -d'_{i} \geq 0$.  So $u'_{\sigma}$ is
in the dual cone $\sigma^{\vee}$.  In particular, $\lambda^{\min}
(u'_{\sigma})$ and $\lambda^{\max} (u'_{\sigma})$ are well-defined.

Finally, with $\Delta = \conv \{ 0, u_{1}, \ldots, u_{s} \}$, note that \[ m
\Delta \, = \, \{ u \in \sigma^{\vee} \ | \ \lambda^{\min}(u) \leq m \}.  \]

The distinction between $\lambda^{\min}$ and $\lambda^{\max}$ is
meaningful only in the nonsimplicial case; when $\sigma$ is
simplicial, then the primitive generators of the rays of
$\sigma^{\vee}$ are linearly independent, so the expression $u =
a_{1}u_{1} + \cdots + a_{n} u_{n}$ is unique.

\vspace{5 pt}

In order to show that $m \Delta$ contains a generating set for
$\sigma^{\vee} \cap M$, one seeks lower bounds for $m$.  To get a
rough idea of how one might get such bounds, imagine that $P_{D}$ is
very large, as it will be under the hypotheses of Theorem 2.  When we
add a small, minus-effective divisor $D'$ to $D$, we get $P_{D+D'}$ by
moving the faces of $P_{D}$ in a small distance.  The main idea is to
control the decrease in the lengths of the edges as the faces move in. 
After the faces containing $u_{\sigma}$ move in a small distance, the
new vertex $u_{\sigma} + u'_{\sigma}$ of $P_{D+D'}$ will be inside
$P_{D}$ and a small distance from the old vertex $u_{\sigma}$ of
$P_{D}$.  We can measure this distance by $\lambda^{\min}
(u'_{\sigma})$.  Suppose that $P_{D}$ contains $u_{\sigma} + t \Delta$
for some large positive $t$.  Looking out from the new vertex
$u_{\sigma} + u'_{\sigma}$ in the direction of the ray spanned by
$u_{i}$, we see that $P_{D}$ contains $u_{\sigma} + u'_{\sigma} + b
u_{i}$ for $0 \leq b \leq t - \lambda^{\min}(u'_{\sigma}).$ Now we
want to know what portion of this segment is actually contained in
$P_{D+D'}$.  This will depend on how far the faces cutting off the
other end of the edge move in.  If these faces move in a distance $r$
with respect to the dual lattice, then the resulting edge of
$P_{D+D'}$ will have length at least $t - \lambda^{\min}(u'_{\sigma})
- r$.  The key to giving lower bounds for $m$ will be the following
proposition, which makes the essence of this discussion precise.

\vspace{10 pt}

\noindent \textbf{Proposition} \emph{ Let $X$ be a complete toric
variety, and $\sigma$ a maximal cone in the fan defining $X$.  Let
$D$ and $D'$ be $T$-\Q -Cartier divisors such that $D$ is nef and $0 \geq D'
\geq K_{X}$.  Let $t = \min \{ D \cdot V(\sigma \cap \tau) \}$ and $m
= \min \{ (D + D') \cdot V(\sigma \cap \tau) \}$, where $\tau$ varies
over all maximal cones adjacent to $\sigma$.  Suppose $t \geq
\lambda^{\min}(u'_{\sigma})$.  Then \[ m \ \geq \ t - \lambda^{\min}
(u'_{\sigma}) - 1.  \] }

Proof. Although we will only use the proposition as stated, we will
prove somewhat more.  We replace the global condition $0 \geq D' \geq
K_{X}$ by the following ``local'' conditions near $\sigma$:
\begin{enumerate}
\item The restriction of $D'$ to $U_{\sigma}$ is minus-effective, 
i.e., $d'_{i} \leq 0$ for each primitive generator $v_{i}$ of a ray of 
$\sigma$. 
\item There is a positive rational number $r$ and, for each
maximal cone $\tau$ adjacent to $\sigma$, a primitive
generator $v_{j}$ of a ray in $\tau \setminus \sigma$ such that
$d'_{j} \geq -r$.
\end{enumerate}
Under these revised hypotheses, we will show that $m \geq t -
\lambda^{\min}(u'_{\sigma}) - r.$ In the case where $0 \geq D' \geq
K_{X}$, the conditions hold for $r = 1$, so the proposition as stated 
will follow.

Let $\tau$ be a maximal cone adjacent to $\sigma$.  Let $u$ be the
primitive generator of the ray of $\sigma^{\vee}$ perpendicular to
$\sigma \cap \tau$, and let $v_{j}$ be the primitive generator of a
ray in $\tau \setminus \sigma$ such that $d'_{j} \geq -r$.  Let $c = t
- \lambda^{\min}(u'_{\sigma})$, and let $k = (D+D') \cdot V(\sigma
\cap \tau)$.  We aim to show $c-k \leq r.$

First, we claim that $u_{\sigma} + u'_{\sigma} + cu$ is in $P_{D}$. 
Since $D$ is nef, $P_{D}$ contains $u_{\sigma} + t \Delta$, so it will
suffice to show $\lambda^{\min}(u'_{\sigma} + cu) \leq t$.  Now, $cu$
is in $\sigma^{\vee}$ and $\lambda^{\min}(cu) = c$, so \[
\lambda^{\min}(u'_{\sigma} + cu) \ \leq \ \lambda^{\min}(u'_{\sigma}) + c
\ = \ t.  \] This proves that $u_{\sigma} + u'_{\sigma} + cu$ is in 
$P_{D}$.  Therefore, we have
\begin{equation}
    \<u_{\sigma} + u'_{\sigma} + cu, v_{j} \> \ \geq \ -d_{j}.
\end{equation}

Next, recall that $u_{\sigma} + u'_{\sigma}$ and $u_{\tau} +
u'_{\tau}$ are the points of $M_{\Q}$ associated to $D+D'$ for
$\sigma$ and $\tau$, respectively.  By Lemma \ref{Laterveer}, \[
u_{\tau} + u'_{\tau} \ = \ u_\sigma + u'_{\sigma} + ku.  \] So,
\begin{equation} \< u_\sigma + u'_{\sigma} + k u , v_{j} \> \ = \ \<
u_{\tau} + u'_{\tau}, v_{j} \> \ = \ -d_{j} -d_{j}' \ \leq \ -d_{j} + r.
\end{equation} 
Subtracting (2) from (1) we have $(c-k) \<u , v_{j} \> \geq -r$. 
Since $\< u,v_{j} \>$ is a negative integer, it follows that $c - k
\leq r$.  \hfill $\Box$

\vspace{5 pt}

Remark. The conclusion of the proposition is false in general if $t
< \lambda^{\min}(u'_{\sigma})$.  Consider, for example, the complete
toric surface $X$ whose fan is spanned by three rays, the primitive
generators of which satisfy $v_{1} + v_{2} + 2v_{3} = 0$. ($X$
is isomorphic to the weighted projective plane $\P(1,1,2)$.)  Let
$\sigma$ be the cone spanned by $v_{2}$ and $v_{3}$.  Taking $T$-\Q
-Cartier divisors $D = D_{1}$ and $D' = K_{X}$, one computes $t = D
\cdot D_{2} = (1/2)$, $\lambda^{\min}(u'_{\sigma}) = 2$, and
$(D+D') \cdot D_{2} = -3$, which is strictly less than $(1/2) -
2 - 1.$

\vspace{5 pt}

The proposition gives good lower bounds for $m$, provided we can
give good upper bounds for $\lambda^{\min}(u'_{\sigma}).$  We will get
sufficient bounds indirectly by using convexity to bound
$\lambda^{\max}(u'_{\sigma})$.

\begin{lemma} \label{convexity bound} Let $D'$ be $T$-\Q -Cartier, with
    $0 \geq D' \geq K_{X}$.  Then $\lambda^{\max}(u'_{\sigma}) \leq
    \lambda^{\max}(u)$ for any lattice point $u$ in the interior of
    $\sigma^{\vee}$.
\end{lemma}

Proof: For any lattice point $u$ in the interior of $\sigma^{\vee}$,
and for the primitive generator $v_{j}$ of any ray of $\sigma$, $\< u,
v_{j} \>$ is a positive integer.  Now $\< u'_{\sigma}, v_{j} \> =
-d'_{j}$, which, since $D' \geq K_{X}$, is at most 1.  Therefore $u -
u'_{\sigma}$ is in $\sigma^{\vee}$.  Since $\lambda^{\max}$ is convex
and nonnegative on $\sigma^{\vee}$, it follows that
$\lambda^{\max}(u'_{\sigma}) \leq \lambda^{\max}(u)$.  \hfill $\Box$

\section{Proof of Theorem 2}
 
Let $X, D,$ and $D'$ satisfy the hypotheses of Theorem 2, and let
$\sigma$ be a maximal cone in the fan defining $X$.  Let $m = \min \{
(D + D') \cdot V(\sigma \cap \tau) \}$, where $\tau$ varies over all
maximal cones adjacent to $\sigma$.  Let $\{u_{1}, \ldots, u_{s}\}$ be
the primitive generators of the rays of $\sigma^{\vee}$, and let
$\Delta = \conv \{ 0, u_{1}, \ldots, u_{s} \}$.  To prove Theorem 2,
it will suffice to show that $m \Delta$ contains a generating set for
$\sigma^{\vee} \cap M$.  To prove this, we will give a canonical
subdivision of $\sigma^{\vee}$ and show that, for each maximal cone
$\gamma$ of the subdivision, $\gamma \cap m \Delta$ contains a
generating set for $\gamma \cap M$.

We claim that $\lambda^{\max}$ is piecewise-linear and therefore
defines a canonical subdivision of $\sigma^{\vee}$: the subdivision
whose maximal cones are the maximal subcones of $\sigma^{\vee}$ on
which $\lambda^{\max}$ is linear.  This subdivision can also be
realized by looking at $Q = \conv \{ u_{1}, \ldots, u_{s} \}$ and
taking the cones over the ``lower faces'' of $Q$, i.e., the faces of
$Q$ visible from the vertex 0 of $\sigma^{\vee}$.  Indeed, for any $t
> 0$, $tQ$ is the set of $u$ in $\sigma^{\vee}$ that can be written $u
= a_{1}u_{1} + \cdots + a_{s} u_{s}$ with $a_{i} \geq 0$ and $a_{1} +
\cdots + a_{s} = t$.  Now the points in the lower faces of $Q$ are
precisely those points that are not contained in $tQ$ for any $t > 1$. 
So the restriction of $\lambda^{\max}$ to the lower faces of $Q$ is
identically 1.  Since $\lambda^{\max}(cu) = c \lambda^{\max}(u)$ for
any $c \geq 0$, it follows that $\lambda^{\max}$ is linear precisely
on the cones over the lower faces of $Q$.

Let $\gamma$ be the cone over a maximal lower face of $Q$, and let
$\gamma^{(1)} \subset \{u_{1}, \ldots, u_{s} \}$ denote the set of
primitive generators of the rays of $\gamma$.  We must show that
$\gamma \cap m \Delta$ contains a generating set for $\gamma \cap M$. 
Every point of $\gamma$ can be written as a nonnegative linear
combination: \[ u \, = \, a_{1} u_{i_{1}} + \cdots + a_{n} u_{i_{n}}, \]
where $a_{j} \geq 0$, and $\{ u_{i_{j}} \} \subset \gamma^{(1)} $ is
linearly independent.  This expression can be decomposed as a
nonnegative integer combination of the $\{ u_{i_{j}} \}$ plus a
nonnegative fractional part.  So $\gamma \cap M$ is generated by $0$
and $\gamma^{(1)}$ together with $\{ (a_{1} u_{i_{1}} + \cdots + a_{n}
u_{i_{n}}) \in M \ | \ 0 \leq a_{j} < 1, \{ u_{i_{j}} \} \subset
\gamma^{(1)}$ linearly independent\}.  By Lemma \ref{cor to Fuj}, $m
\Delta$ contains $0$ and $\gamma^{(1)}$.  It will therefore suffice to
show that any lattice point $p$ that is a nonnegative fractional linear
combination of some independent set $\{ u_{i_{j}} \} \subset
\gamma^{(1)}$ is contained in $m \Delta$.  For this, it will suffice 
to show that $m \geq \lambda^{\max}(p)$.

Suppose $p = a_{1}u_{i_{1}} + \cdots + a_{n} u_{i_{n}} \in M$, where
$0 \leq a_{j} < 1$, and $\{ u_{i_{j}} \} \subset \gamma^{(1)}$ is
linearly independent.  Then $p' = (1-a_{1})u_{i_{1}} + \cdots +
(1-a_{n})u_{i_{n}}$ is a lattice point in the interior of
$\sigma^{\vee}$.  By Lemma \ref{convexity bound}, $\lambda^{\max}
(u'_{\sigma}) \leq \lambda^{\max} (p'),$ and since $\lambda^{\max}$ is
linear on $\gamma$ and $\lambda^{\max}(u_{i_{j}}) = 1$, we have \[
\lambda^{\max} (p') \ = \ (1-a_{1}) + \cdots + (1-a_{n}) \ = \ n -
\lambda^{\max} (p).  \] Therefore, \[ \lambda^{\max} (u'_{\sigma}) \
\leq \ n - \lambda^{\max} (p).  \] Let $t$ be as in the Proposition,
i.e., $t = \min \{ D \cdot V(\sigma \cap \tau) \}$, where $\tau$ varies
over all maximal cones adjacent to $\sigma$.  Then $t \geq n + 1 >
\lambda^{\min} (u'_{\sigma})$, so we can apply the Proposition with $r
= 1$ to obtain
\begin{eqnarray*}
  m & \geq & n + 1 - \lambda^{\min} (u'_{\sigma}) - 1. \\
    & \geq & n - \lambda^{\max} (u'_{\sigma}). \\
    & \geq & \lambda^{\max} (p).
\end{eqnarray*}

\vspace{-25 pt} \hfill $\Box$

\vspace{25 pt}

Remark. The collection of cones over the lower faces of $Q$ is an
example of what is called a ``regular subdivision''.  In general, a
regular subdivision of a cone is constructed by choosing a nonzero
point on each of the rays of the cone, and perhaps specifying some
additional rays inside the cone with nonzero points on them as well. 
One looks at the convex hull of all of these points and then takes the
cones over all of the lower faces.  For more details on regular
subdivisions of convex polytopes, see \cite{Lee} and \cite{Ziegler}. 
The translation from polytopes to cones is straightforward.

In the toric literature, regular subdivisions have generally been
applied to the fan defining a toric variety, and sometimes to the
polytope defining an ample line bundle.  See, for instance, \cite{Oda
and Park}, \cite[Chapter 7]{GKZ}, and \cite[\S I.2]{Toroidal 
Embeddings}.  The regular subdivisions that we have used
in this paper are of the dual cones $\{ \sigma^{\vee} \}$.  The author
is not aware of any significant geometric interpretation for these
subdivisions.

The subdivisions of a fan $\Sigma$ correspond naturally and
bijectively to the proper birational toric morphisms $\tilde{X}
\rightarrow X(\Sigma)$ \cite[Section 2.5]{Fulton}, and the regular
subdivisions of $\Sigma$ are precisely those for which the
corresponding morphism is projective.  A regular subdivision of
$\Sigma$ is obtained by specifying a continuous function $\Psi$ on the
support of $\Sigma$ that is convex and piecewise-linear on each cone. 
By subdividing each cone of $\Sigma$ into the maximal subcones on
which $\Psi$ is linear, we get a projective birational morphism for
which $\Psi$ is the piecewise-linear function associated to a
relatively ample $T$-\Q -Cartier divisor on $\tilde{X}$.  In
particular, if $D = \sum d_{i}D_{i}$ is a $T$-\Q -Weil divisor on $X$,
and if we define $\Psi^{\max}_{D}$ on each maximal cone $\sigma$ by \[
\Psi^{\max}_{D}(v) = \max \left\{ \sum_{v_{i} \in \sigma} -a_{i}d_{i}
\ | \ \sum_{v_{i} \in \sigma} a_{i}v_{i} = v, a_{i} \geq 0 \right\}, \]
then we get the unique projective birational morphism $\pi: \tilde{X}
\rightarrow X$ such that the proper transform of $D$ is \Q -Cartier
and relatively ample, and $\pi$ is an isomorphism in codimension 1. 
If $D$ is effective (resp.\ minus effective) then, for each maximal
cone $\sigma$, the same subdivision is obtained by looking at $\conv
\{ (1/d_{i}) v_{i} \ | \ v_{i} \in \sigma \}$ (resp.  $\conv \{
-(1/d_{i}) v_{i} \ | \ v_{i} \in \sigma \} $) and taking the cones
over the upper faces (resp.  lower faces).  Note that, for a
subdivision of a fan to be regular, it is not enough for the
subdivision to be regular on each cone.  This is the toric
manifestation of the fact that quasiprojectivity is not local on the
base (see \cite[II.5.3]{EGA}).

\section{Appendix: Proof of Fujino's Theorem$^{+}$}

The ideas and techniques of the main part of this paper also give a
new proof of Fujino's Theorem$^{+}$.  This yields a unified
combinatorial approach to Fujita's conjectures for toric varieties
with arbitrary singularities, which is independent of vanishing theorems and
toric Mori theory.

From the Proposition, we can immediately deduce a generalization of
Fujita's freeness conjecture for toric varieties with arbitrary
singularities.  The statement we get in this way is similar to
Fujino's Theorem$^{+}$, but without the improved bound obtained by
excluding the extremal case where $X$ is $\P^{n}$.

\vspace{10 pt}

\noindent \textbf{Corollary to Proposition} \emph{ Let $X$ be a projective 
$n$-dimensional toric variety.  Let $D$ and $D'$ be \Q -Cartier divisors 
such that $0 \geq D' \geq K_{X}$ and $D \cdot C \geq n+1$ for all 
$T$-curves $C$.  Then $D +D'$ is nef.}

\vspace{10 pt}

Proof. We may assume that $D$ and $D'$ are $T$-\Q -Cartier.  By the
Proposition, it suffices to show that $\lambda^{\min}(u'_{\sigma})
\leq n$.  Write $u'_{\sigma} = a_{1} u_{1} + \cdots + a_{n}u_{n}$,
where $a_{i} \geq 0$ and the $u_{i}$ are linearly independent
primitive generators of rays of $\sigma^{\vee}$.  The condition $D'
\geq K_{X}$ implies that each $a_{i} \leq 1$.  So
$\lambda^{\min}(u'_{\sigma}) \leq a_{1} + \cdots + a_{n} \leq n$. 

\nopagebreak

\hfill $\Box$

\vspace{5 pt}

To prove Fujino's Theorem$^{+}$, it remains to show that the
bound on the intersection numbers can be improved by one by excluding
the case where $X \cong \P^{n}$.

Using the Proposition, we can work ``locally,'' considering one 
maximal cone $\sigma$ at a time.  The following lemma allows us 
to reduce to the case where $\sigma$ is regular.

\begin{lemma} \label{not regular} Let $D'$ be a $T$-\Q -Cartier
divisor, $0 \geq D' \geq K_{X}$.  If $\sigma$ is a maximal cone in the
fan defining $X$ that is not regular, then $\lambda^{\min}
(u'_{\sigma}) \leq n-1$.
\end{lemma}

Proof. By Lemma \ref{convexity bound}, it will suffice to show that
there is a lattice point $p$ in the interior of $\sigma^{\vee}$ such
that $\lambda^{\max}(p) \leq n-1$.  Since $\sigma$ is not regular,
$\sigma^{\vee}$ is not regular either.  Consider two cases, according
to whether $\sigma^{\vee}$ is simplicial.

Suppose $\sigma^{\vee}$ is simplicial, and let $u_{1}, \ldots, u_{n}$
be the primitive generators of the rays of $\sigma^{\vee}$.  Since
$\sigma^{\vee}$ is not regular, after possibly renumbering the
$u_{i}$, there is a lattice point $u = a_{1}u_{1} + \cdots +
a_{r}u_{r}$ in $M$, where $0 < a_{i} < 1$ and $r \geq 2$.  Then $p =
u_{1} + \cdots + u_{n} - u$ and $p' = u + u_{r+1} + \cdots + u_{n}$ are
lattice points in the interior of $\sigma^{\vee}$.  Now
$\lambda^{\max}(p) + \lambda^{\max}(p') = 2n - r \leq 2n -2$.  So
$\min \{ \lambda^{\max}(p), \lambda^{\max}(p') \} \leq n - 1$, as
required.

Suppose $\sigma^{\vee}$ is not simplicial.  Let $\{u_{1}, \ldots,
u_{n} \}$ be a set of linearly independent primitive generators of
rays in some subcone of $\sigma^{\vee}$ on which $\lambda^{\max}$ is
linear.  Let $\gamma$ be the cone spanned by $\{ u_{1}, \ldots, u_{n}
\}$.  Since $\gamma \subsetneq \sigma^{\vee}$, at least one of the
facets of $\gamma$ is not contained in a face of $\sigma^{\vee}$.  Say
$\tau$, spanned by $\{ u_{1}, \ldots, u_{n-1} \}$, is not contained in
a face of $\sigma^{\vee}$.  Then the relative interior of $\tau$ is
contained in the interior of $\sigma^{\vee}$.  In particular, $p =
u_{1} + \cdots + u_{n-1}$ is a lattice point in the interior of
$\sigma^{\vee}$.  Since $\lambda^{\max}$ is linear on $\gamma$ and
$\lambda^{\max}(u_{i}) = 1$, it follows that $\lambda^{\max}(p) =
n-1$.  \hfill $\Box$

\vspace{10 pt}

\noindent \textbf{Proof of Fujino's Theorem$^{+}$.} Let $X, D$, and
$D'$ satisfy the hypotheses of Fujino's Theorem$^{+}$.  Let $\sigma$
be a maximal cone in the fan defining $X$.  Let $t = \min \{ D \cdot
V(\sigma \cap \tau) \}$ and $m = \min \{ (D+D') \cdot V(\sigma \cap
\tau) \}$, where $\tau$ varies over all maximal cones adjacent to
$\sigma$.  It will suffice to show that $m \geq 0$.  If $\sigma$ is
not regular, then, by Lemma \ref{not regular}, $\lambda^{\min}
(u'_{\sigma}) \leq n-1$.  Applying the Proposition, we have $ m \geq t
- n \geq 0$, as required.

We may therefore assume that $\sigma$ is regular.  Let $v_{1}, \ldots,
v_{n}$ be the primitive generators of the rays of $\sigma$, and let
$u_{1}, \ldots, u_{n}$ be the dual basis.  In particular, $u_{1},
\ldots, u_{n}$ are the primitive generators of the rays of
$\sigma^{\vee}$.  By adding to $D$ the numerically trivial divisor
$\sum \< u_{\sigma} , v_{i} \> D_{i}$, we may assume that $u_{\sigma}
= 0$, and hence $P_{D}^{\sigma} = P_{D}$.  Consider two cases,
according to whether $P_{D}$ is a simplex.

\vspace{5 pt}

\emph{Case {\rm 1:} $P_{D}$ is a simplex.} In this case, $\Pic(X) \cong \Z$. 
We may linearly order $\Pic X \otimes \Q$ so that $[D]$ is positive
and the nef divisor classes are exactly those that are greater than
or equal to zero.  We aim to show that $[D + D'] \geq 0$.  It will
suffice to show that $P_{D+D'}$ is nonempty.  In fact, we will show
that $p = u_{1} + \cdots + u_{n}$ is in $P_{D+D'}$.  Furthermore,
since $[D] \geq [(n/t)D]$, it will suffice to prove this in the
case where $t = n$.

Let $v_{0}$ be the primitive generator of the unique ray of the fan
defining $X$ that is not in $\sigma$.  Now $P_{D+D'}$ is cut out by
the inequalities $\< u, v_{i} \> \geq -d'_{i}$ for $1 \leq i \leq n$,
and $\< u, v_{0} \> \geq -d_{0} -d'_{0}$.  For $1 \leq i \leq n$, we
have $\< p, v_{i} \> = 1 \geq -d'_{i}$.  So it will suffice to show
that $\< p, v_{0} \> \geq -d_{0} - d'_{0}$.

Write $P_{D} = \conv \{ 0, a_{1} u_{1}, \ldots, a_{n}u_{n} \}$.  After
possibly renumbering, we may assume $a_{1} = \min \{ a_{i} \} = n$. 
Furthermore, one of the $a_{i}$ must be strictly greater than $n$
(otherwise $P_{D}$ would be a regular simplex and so $X$ would be
isomorphic to $\P^{n}$).  The ray spanned by $v_{0}$, which is
perpendicular to the face of $P_{D}$ not containing 0, is also spanned
by \[ v \ = \ -a_{2} \cdots a_{n} v_{1} - \cdots - a_{1} \cdots
\hat{a_{i}} \cdots a_{n} v_{i} - \cdots - a_{1} \cdots a_{n-1}
v_{n}.\] So $ v_{0} = bv $ for some positive rational number $b$. 
Since $v_{0}$ is a lattice point, $b a_{1} \cdots \hat{a_{i}} \cdots
a_{n}$ must be an integer for each $i$.  In particular, $ba_{2} \cdots 
a_{n}$ is an integer.  Therefore, \[ d_{0} \ = \ - \<
a_{1} u_{1}, v_{0} \> \ = \ b a_{1} \cdots a_{n} \ = \ n b a_{2} \cdots
a_{n}  \] is an integer.

Now, since $a_{1} = \min \{ a_{i} \} = n$ and some $a_{i} > n$, we
have \[ \< p, v_{0} \> \ = \ b \< p, v \> \ > \ - n b a_{2} \cdots
a_{n} \ = \ - d_{0}.  \] Since both $\< p, v_{0} \>$ and $- d_{0}$ are
integers, their difference must be at least 1.  So $\< p, v_{0} \>
\geq -d_{0} + 1 \geq -d_{0} - d'_{0}$.

\vspace{5 pt}

\emph{Case {\rm 2:} $P_{D}$ is not a simplex.} Since $\sigma$ is simplicial
but $P_{D}$ is not a simplex, $P_{D}$ has a vertex $u_{0}$ that is not
adjacent to $u_{\sigma} = 0$.  Define a piecewise linear function
$\lambda$ on $\sigma^{\vee}$ by \[ \lambda(u) = \min \{ (na_{0} +
a_{1} + \cdots + a_{n}) \ | \ a_{0} u_{0} + \cdots + a_{n} u_{n} = u,
a_{i} \geq 0 \}.  \] Now $P_{D}$ contains $\conv \{0, u_{0}, nu_{1},
\ldots, nu_{n} \} = \{ u \in \sigma^{\vee} \ | \ \lambda(u) \leq n\}$.  An
argument identical to the proof of the Proposition shows that $m \geq
n - \lambda (u'_{\sigma}) - 1.$ It will therefore suffice to show that
$\lambda (u'_{\sigma}) \leq n - 1$.

After possibly renumbering, we may write $u_{0} = b_{1} u_{1} + \cdots
+ b_{r} u_{r}$, where $b_{i} > 0$ and $r \geq 2$.  We claim that
$b_{i} \geq n$.  Indeed, the rays along the edges of $P_{D}$ coming
out from $u_{0}$ span a translated cone containing $P_{D}$, and hence
containing $0$.  Since $\< u_{0}, v_{i} \> = b_{i} > 0$, there must be
some vertex $u$ of $P_{D}$ adjacent to $u_{0}$ such that $\< u, v_{i}
\> < \< u_{0}, v_{i} \>$.  Let $C$ be the $T$-curve corresponding to
the edge connecting $u$ and $u_{0}$.  Since $u \in P_{D} \subset
\sigma^{\vee}$, we have $b_{i} = \< u_{0}, v_{i} \> \geq \< u_{0} - u,
v_{i} \>$.  Let $u'$ be the primitive generator of the ray spanned by
$u_{0} - u$.  By Lemma \ref{Laterveer}, we have \[ \< u_{0} - u, v_{i}
\> \ \geq \ (D \cdot C) \< u', v_{i} \> \ \geq \ D \cdot C \ \geq \ n. 
\] This proves the claim.

Let $c_{0}$ be the largest rational number such that $u'_{\sigma} -
c_{0}u_{0}$ is in $\sigma^{\vee}$, i.e.,
\[
c_{0} = \min \{ -d'_{i}/b_{i} \ | \ 1 \leq i \leq r \}.
\]
So we may write \[ u'_{\sigma} \ = \ c_{0} u_{0} + \cdots + c_{n}
u_{n}, \] where $c_{i} \geq 0$ and some $c_{i} = 0$ for $i \geq 1$.
Say $c_{1} = 0$.  We claim that $nc_{0} + c_{2} \leq 1$.  Indeed, \[
nc_{0} + c_{2} \ \leq \ b_{2} c_{0} + c_{2} \ = \ -d'_{2} \ \leq \
1.\] Therefore,
\[
    \lambda(u'_{\sigma})  \ \leq \ nc_{0} + c_{2} + \cdots + c_{n}
                          \ \leq \ 1 - d'_{3} - \cdots - d'_{n}
			  \ \leq \ n - 1.
\]

\vspace{-25 pt} \hfill $\Box$

\begin{flushleft}
\textsc{Department of Mathematics\\ University of Michigan\\ 2074 East
Hall\\ Ann Arbor, MI 48109\\ USA}

\vspace{10 pt}

\textit{E-mail address: sdpayne@umich.edu}
\end{flushleft}

\end{document}